\documentclass[graybox]{svmult}

\usepackage{mathptmx}       
\usepackage{helvet}         
\usepackage{courier}        
\usepackage{type1cm}        
\usepackage{makeidx}         
\usepackage{graphicx}        
\usepackage{multicol}        
\usepackage[bottom]{footmisc}

\usepackage{amsfonts}
\usepackage{amsmath}
\usepackage{amssymb}
\usepackage{graphicx}
\usepackage{enumerate}
\usepackage{multicol}
\usepackage{mathrsfs}
\usepackage[usenames,dvipsnames]{pstricks}
\usepackage{epsfig}
\usepackage{pst-grad} 
\usepackage{pst-plot} 

\makeindex             

\begin{document}

\title*{The rolling problem: overview and challenges}
\author{Yacine Chitour, Mauricio Godoy Molina and Petri Kokkonen}
\institute{Yacine Chitour \at L2S, Universit\'e Paris-Sud XI, CNRS and Sup\'elec, Gif-sur-Yvette, 91192, France, \email{yacine.chitour@lss.supelec.fr}
\and Mauricio Godoy Molina \at L2S, Universit\'e Paris-Sud XI, CNRS and Sup\'elec, Gif-sur-Yvette, 91192, France, \email{mauricio.godoy@gmail.com}
\and Petri Kokkonen \at L2S, Universit\'e Paris-Sud XI, CNRS and Sup\'elec, Gif-sur-Yvette, 91192, France and University of Eastern Finland, Department of Applied Physics, 70211, Kuopio, Finland,
\email{petri.kokkonen@lss.supelec.fr}}
%
%
\maketitle

\abstract*{In the present paper we give a historical account --ranging from classical to modern results-- of the problem of rolling two Riemannian manifolds one on the other, with the restrictions that they cannot instantaneously slip or spin one with respect to the other. On the way we show how this problem has profited from the development of intrinsic Riemannian geometry, from geometric control theory and sub-Riemannian geometry. We also mention how other areas --such as robotics and interpolation theory-- have employed the rolling model.}

\abstract{In the present paper we give a historical account --ranging from classical to modern results-- of the problem of rolling two Riemannian manifolds one on the other, with the restrictions that they cannot instantaneously slip or spin one with respect to the other. On the way we show how this problem has profited from the development of intrinsic Riemannian geometry, from geometric control theory and sub-Riemannian geometry. We also mention how other areas --such as robotics and interpolation theory-- have employed the rolling model.}


\section{Introduction}\label{sec:intro}

Differential geometry has been inextricably related to classical mechanics, since its very conception in the 18th century. As a matter of fact, back in the days, this area of research was referred to as rational mechanics. 
The basic idea of this point of view is reasonably simple: to a given mechanical system ${\mathbb M}$, one can associate a differentiable manifold $M$ in such a way that each possible state of the system corresponds to a unique point in $M$. In this way, each possible velocity vector of ${\mathbb M}$ at a given configuration is represented as a tangent vector to $M$ at the corresponding point. The classical dictionary goes as follows:
\begin{enumerate}
\item Physical data (such as masses, lengths, etc.) of elements in ${\mathbb M}$ induce a Riemannian metric in $M$ representing the kinetic energy.
\item Linear restrictions imposed on the positions of ${\mathbb M}$ (or that can be integrated to such) translate to submanifolds of $M$.
\end{enumerate}

In the late 19th century, physicists noted there were plenty of mechanical systems not considered by the above dictionary. These systems were named non-holonomic, opposed to holonomic systems which are defined in the second point of the dictionary above. A mechanical system ${\mathbb M}$ is non-holonomic if its dynamics has linear restrictions that cannot be integrated to constraints of the position. For various examples and a brief historical bibliography, we refer the interested reader to the survey~\cite{bloch05}. A well-known early example of these systems is the sphere rolling on the plane without sliding or spinning, studied (with some variants) by S. A. Chaplygin in the seminal works~\cite{Chap1,Chap2}. Our aim in this paper is to give a general look at some of the most important breakthroughs in mathematics that gave us some understanding of the generalized version of this system consisting on two Riemannian manifolds $M$ and $\hat M$ of the same dimension rolling one against the other, not allowing instantaneous spins or slips. Nowadays these systems are often studied in connection to sub-Riemannian and Riemannian geometry~\cite{montgomery06,sharpe97} and geometric control theory~\cite{agrachev04}.

The structure of the paper is the following. In Section~\ref{sec:old} we recall two major players in the study of the mechanical system described above and early differential geometry: S. A. Chaplygin and \'E. Cartan. Chaplygin studies for the first time the problem from a mechanical point of view and finds first integrals of motion in different situations. Cartan's development and his celebrated ``five variables'' paper were not evidently connected to the rolling model at the time of their publication, see~\cite{bryant06}, nevertheless we present them from our point of view. In Section~\ref{sec:nomizu}, we briefly present Nomizu's breakthrough introduction of the dynamics of rolling in higher dimensions, through embedded submanifolds of Euclidean space and its relation to Cartan's development. In Section~\ref{sec:2dim} we present how the problem was brought back to life when control theory sees in differential geometry a useful tool to treat the controllability issue of the rolling model in two dimensions and some geometric consequences of optimality conditions. Section~\ref{sec:higerdim} surveys how the higher dimensional rolling system was re-discovered and how it appears naturally in geometric interpolation. Finally in Section~\ref{sec:intrinsic} we present the latest results that have been obtained concerning the controllability of the system and its symmetries. We conclude with a brief discussion on some generalizations and open problems.

\section{The early years: Mechanics and the new differential geometry}\label{sec:old}

The first time the problem of a ball rolling on the plane was considered as worthy of study was in the seminal papers of S. A. Chaplygin~\cite{Chap1,Chap2}, one of the fathers of non-holonomic mechanics. The results were considered surprisingly difficult at the time, and for~\cite{Chap1} Chaplygin won the Gold Medal of Russian Academy of Sciences. The main results he obtained were first integrals of motion for the system in several geometric situations. Even these seemingly elementary problems contains unexpected difficulties and bottlenecks when trying to obtain closed formulae for the dynamics. As stated in~\cite{Chap1}, after observing that the differential equation of the dynamics can be integrated in quadratures. Essentially at the same time, \'E. Cartan was developing his coordinate-free differential geometry. With this new language he was able to propose and study many problems, most often related to the search of invariants of geometric systems. In this survey, we will only focus in two of his many ideas: the search for invariants and symmetries for control systems with two controls and five degrees of freedom, and the definition of affine Riemannian holonomy through the development of a curve. Both of this ideas will appear several other times in this survey.

\subsection{Chaplygin's ball}\label{subsec:chaplygin}

In the year 1897 the work~\cite{Chap1} written by S. A. Chaplygin was published. This papers is one of a series of research articles in which Chaplygin analyzed non-holonomic systems. Also of particular relevance to this survey is another paper~\cite{Chap2}. In particular he was interested in studying first integrals and equations of motion for different systems of rolling balls.

To illustrate his results, Chaplygin was able to find an integral of motion for the system of a homogeneous small ball of mass $m_1$ and a homogeneous sphere of mass $m_2$, in which the ball rolls without slipping inside the sphere. We will think of the dynamics occurring in Euclidean 3-space. Let $O$ be the center of the sphere, let $G$ be the center of the moving ball and $A$ the point of contact between the two. Introducing the quantities $a={\rm dist}(O,G)$ and $b={\rm dist}(O,A)$, then one has the integrals of motion:
\[
\sum_{i=1}^2 m_i\left(y_i\frac{dx_i}{dt}-x_i\frac{dy_i}{dt}\right)+M\left(\frac{b}{a}-1\right)\left(\beta\frac{d\alpha}{dt}-\alpha\frac{d\beta}{dt}\right)={\rm const.}
\]
Where $A=(\alpha,\beta,\gamma)$ with respect to a fixed frame $OX'Y'Z'$, and the points $G=(x_1,y_1,z_1)$ and $O=(x_2,y_2,z_2)$ with respect to a moving frame $AXYZ$, with axes at all times parallel to those in $OX'Y'Z'$. Additionally $M=m_1+m_2$ denotes the mass of the system. 

The equations of motion are complicated and it serves little purpose to write them down here. Nevertheless, there is an interesting historical remark at this point. After arriving at a very complicated differential equation to describe the dynamics of the system, Chaplygin observes it can be written in the form
\[
\frac{dv}{d\zeta}+v\Phi(\zeta)+\Psi(\zeta)=0,
\]
for some appropriate functions $\Phi$ and $\Psi$ after a series of changes of variables. He then ventures to say
\begin{quote}
[\ldots] and, therefore, can be integrated in quadratures. We will not write out these quadratures since they are rather cumbersome.
\end{quote}
As far as we know, the integration of differential equations connected to the problem of rolling balls is still an area of active research, see for example~\cite{borisov12}.

\subsection{Cartan's ``five variables'' paper}\label{subsec:1/3}

A rank $l$ vector distribution $D$ on an $n$-dimensional manifold $M$ or $(l, n)$-distribution (where $l < n$) is, by definition, an $l$-dimensional subbundle of the tangent bundle $TM$, i.e., a smooth assignment $q\mapsto D|_q$ defined on $M$ where $D|_q$ is an $l$-dimensional subspace of the tangent space $T_qM$. Two vector distributions $D_1$ and $D_2$ are said to be equivalent, if there exists a diffeomorphism $F : M \rightarrow M$ such that $F_*D_1|_q = D_2|_{F(q)}$ for every $q\in M$. Local equivalence of two distributions is defined analogously.

Cartan's equivalence problem consists in constructing invariants of distributions with respect to the equivalence relation defined above. 
A seminal contribution by \'E. Cartan in \cite{cartan10}  was the introduction of the ``reduction-prolongation'' procedure for building invariants and the characterization for $(2,5)$-distributions via a functional invariant (Cartan's tensor) which vanishes precisely when the distribution is flat, that is, when it is locally equivalent to the (unique) graded nilpotent Lie algebra ${\mathfrak h}$ of step $3$ with growth vector $(2,3,5)$.

In the same paper, Cartan also proved that in this system there is hidden a realization of the $14$-dimensional exceptional Lie algebra ${\mathfrak{g}}_2$. To explain where does it appear, let us recall that an infinitesimal symmetry of an $(l,n)-$distribution $D$ is a vector field  $X\in{\rm VF}(M)$ such that $[X,D]\subseteq D$. Now consider the (unique) connected and simply connected nilpotent Lie group $H$ with Lie algebra ${\mathfrak h}$. The two dimensional subspace of ${\mathfrak h}$ that Lie generates it, can be seen as a $(2,5)-$distribution on $H$. In general, a $(2,5)-$distribution that is bracket generating is nowadays known as a Cartan distribution. In this setting, the following theorem takes place.

\begin{theorem}[Cartan 1910]
The Lie algebra of symmetries of the flat Cartan distribution is precisely ${\mathfrak{g}}_2$, and this situation is maximal, that is, for general Cartan distributions the dimension of the Lie algebra of symmetries is $\leq14$.
\end{theorem}

Moreover, Cartan  gave a geometric description of the flat $G_2$-structure as the differential system that describes space curves of constant torsion $2$ or $1/2$ in the standard unit $3$-sphere (see Section 53 in Paragraph XI in \cite{cartan10}.) 

The connection between this studies by Cartan and the rolling problem comes from the fact that the flat situation described above occurs in the problem of two $2$-dimensional spheres rolling one against the other without slipping or spinning, assuming that the ratio of their radii is $1\colon3$, see~\cite{bor06} for some historical notes and a thorough attempt of an explanation for this ratio. In fact, whenever the ratio of their radii is different from $1\colon 3$, the Lie algebra of symmetries becomes ${\mathfrak{so}}(3)\times{\mathfrak{so}}(3)$, thus dropping its dimension to 6. A complete answer to this strange phenomenon as well as a geometric reason for Cartan's tensor was finally given in two remarkable papers~\cite{zelenko061,zelenko062} (cf. also \cite{agrachev07}), where a geometric method for construction of functional invariants of generic germs of $(2, n)$-distribution for arbitrary $n\geq5$ is developed. It has been recently observed in~\cite{nurowski12} that the Lie algebra of symmetries of a system of rolling surfaces can be $\mathfrak{g}_2$ in the case of non-constant Gaussian curvature.

\subsection{Cartan's development}\label{subsec:development}

\'E. Cartan in~\cite{cartan25} defined a geometric operation, that he called development of a manifold onto a tangent space, in order to define holonomy in terms of ``Euclidean displacements'', i.e.,  elements of ${\rm E}(n)$. In his own words:
\begin{quotation}
Quand on d\'eveloppe l'espace de Riemann sur l'espace euclidien tangent en $A$ le long d'un cycle partant de $A$ et y revenant, cet espace euclidien subit un d\'eplacement et tous les d\'eplacements correspondant aux diff\'erents cycles possibles forment un groupe, appel\'e groupe d'holonomie.
\end{quotation}

An interpretation of this quote in terms of manifolds rolling follows naturally. For a given loop $\gamma\colon[0,\tau]\to M$ on an $n$ dimensional Riemannian manifold $M$, one can roll $M$ against the Euclidean space ${\mathbb R}^n$ obtaining a new curve $\hat\gamma\colon[0,\tau]\to{\mathbb R}^n$, called the development of $\gamma$. By parallel transporting along $\gamma$ any orthonormal frame of $T|_{\gamma(0)}M$, we obtain a rotation $R_\gamma\in{\rm O}(n)$. The fact that $\hat\gamma$ is not necessarily a loop induces a translation $T_\gamma$ corresponding to the vector $\hat\gamma(\tau)-\hat\gamma(0)$. We conclude that we can associate to $\gamma$ an element $(R_\gamma,T_\gamma)$ of the Euclidean group of motions ${\rm E}(n)$. The subgroup ${\rm Hol}^{aff}(M)$ of ${\rm E}(n)$ consisting of all such $(R_\gamma,T_\gamma)$ obtained by rolling along all absolutely continuous loops $\gamma$ is known as the affine holonomy group of $M$ and the orthogonal part ${\rm Hol}(M)\subseteq{\rm O}(n)$ of it is the holonomy group of $M$.

It is known that if $M$ is complete and with irreducible Riemannian holonomy group, the affine holonomy group contains all translations of $T|_xM$, see~\cite[Corollary 7.4, Chapter IV]{kobayashi63}. In other words, under the irreducibility hypothesis, the rotational part of the affine holonomy permits to recover the translational part, and this consists of all the possible translations in $T|_xM$.

Perhaps something that might have been not expected by Cartan is that this concept of development would play a fundamental role in the definition of Brownian motion on a manifold, and the subsequent explosion of interest that stochastic analysis in Riemannian manifolds has had in later decades, see~\cite{Hsu}. For a long time, mathematicians have had the intuition that by rolling an $n$-dimensional manifold $M$ along a given curve $y(t)$ in ${\mathbb R}^n$ with the Euclidean structure, one would obtain a curve in $M$ which resembles the original curve $y(t)$, see~\cite{GG}. The main outstanding idea (as far as we know due to Malliavin) was to use Cartan's development through the orthonormal frame bundle and Wiener's measure, see~\cite{stroock}.

The idea of how to define Brownian trajectories on manifolds is similar to the interpretation given above. Intuitively, one can draw a Brownian path $B(t)$ in ${\mathbb R}^n$, and then one can consider the system of $M$ rolling against ${\mathbb R}^n$ following the path $B(t)$. The precise definition uses a less regular version of Cartan's development and parallel transport.

This naive notion allows one to recover the Laplace-Beltrami operator $\Delta_M$ of the manifold. It is often interpreted as if Brownian paths are the ``integral curves'' for $\Delta_M$. Of course this assertion lacks of mathematical precision, but it introduces the idea that second order differential operators induce ``diffusions'' on the manifold. This point of view has been exploited significantly in the study of stochastic differential equations on manifolds, see~\cite{bismut}.

\section{A ``forgotten'' breakthrough}\label{sec:nomizu}

An important contribution to the understanding of the problem of rolling without slips or spins came to light in the paper~\cite{nomizu78} by K. Nomizu. His aim was to give a mechanical interpretation of certain differential geometric invariants using this system. He mainly focuses in submanifolds of ${\mathbb R}^N$ with the usual Euclidean structure, and so will we along this section.

He begins with a simple general consideration: as a motion occurring in a Euclidean space ${\mathbb R}^N$ without deforming objects, a rolling can be seen as a curve in the Euclidean group ${\rm E}(N)$, that is a function $[0,\tau]\ni t\mapsto f_t\in{\rm E}(N)$ given by
\begin{equation}\label{eq:motion}
f_t=\begin{pmatrix}
C_t&c_t\\0&1
\end{pmatrix},
\end{equation}
where $f_0={\rm Id}$ is the identity matrix of $(N+1)\times(N+1)$, $C_t\in{\rm O}(N)$ and $c_t\in{\mathbb R}^N$. He calls such types of curves $1$-parametric motions.

For a given $1$-parametric motion $\{f_t\}$, he observed that there is a natural time-dependent vector field $X_t$ associated to it. For an arbitrary point $y\in{\mathbb R}^N$ we define $(X_t)_y:=\left.\frac{df_u(x)}{du}\right|_{u=t}$, where $x=f_t^{-1}(y)$. Using equation~\eqref{eq:motion}, one can see that $(X_t)_y=S_ty+v_t$, where $S_t=\frac{dC_t}{dt}C_t^{-1}\in\mathfrak{o}(N)$ and $v_t=-S_tc_t+\frac{dc_t}{dt}\in{\mathbb R}^N$. The corresponding element of the Lie algebra $\mathfrak{e}(N)$
\begin{equation}\label{eq:instmot}
\frac{df_t}{dt}f_t^{-1}=\begin{pmatrix}
S_t&v_t\\
0&0
\end{pmatrix}
\end{equation}
is called the instantaneous motion. Slips and spins can now be encoded in terms of the vector field $X_t$ and the instantaneous motion.
\begin{definition}
The instantaneous motion~\eqref{eq:instmot} is called an instantaneous:
\begin{itemize}
\item standstill if $S_t=0$ and $v_t=0$,
\item translation if $S_t=0$ and $v_t\neq0$,
\item rotation if there exists a point $y_0\in{\mathbb R}^N$ such that $(X_t)_{y_0}=0$ and $S_t\neq0$.
\end{itemize}
\end{definition}

With this at hand, it is possible to define rolling without slipping ({\em skidding} in Nomizu's terminology) nor spinning between $M^n,\hat M^n\hookrightarrow{\mathbb R}^N$.
\begin{definition}
Let $\{f_t\}$ be a $1$-parametric motion such that $f_t(M)$ is tangent to $\hat M$ at a point $y_t\in\hat M$. Assume that $(X_t)_{y_t}=0$ and $S_t\neq0$. The motion $f_t$ is a rolling if for any pair of tangent vectors $X,Y\in T_{y_t}N$
\begin{equation}\label{eq:notwist1}
\langle S_t(X),Y\rangle=0,
\end{equation}
and for any pair of normal vectors $U,V\in T_{y_t}^\bot \hat M$
\begin{equation}\label{eq:notwist2}
\langle S_t(U),V\rangle=0.
\end{equation}
\end{definition}

An equivalent way of stating conditions~\eqref{eq:notwist1} and \eqref{eq:notwist2} is that $S_t$ maps $T_{y_t}\hat M$ to $T_{y_t}^\bot \hat M$ and also maps $T_{y_t}^\bot \hat M$ to $T_{y_t}\hat M$.

This definition allowed Nomizu to find a very concrete realization of Cartan's development. For the case of surfaces rolling on the plane, his result reads
\begin{theorem}[Nomizu 1978]
Let $x_t$ be a smooth curve on a surface $M$ which does not go through a flat point of $M$. There exists a unique rolling $\{f_t\}$ of $M$ on the tangent plane $\Sigma$ at $x_0$ such that $y_t=f_t(x_t)$ is the locus of points of contact. The curve $y_t$ is the development of the curve $x_t$ into $\Sigma$.
\end{theorem}
\vspace{0.5cm}
As a consequence of this result, Nomizu noticed that there is a natural kinematic interpretation of the Levi-Civita connection for a surface $M$, coming from the rolling formulation: a vector field $U(t)$ along the curve $x_t$ is parallel with respect to the Levi-Civita connection of $M$ if and only if $C_t(U(t))$ is a constant vector for all $t$.

As a matter of fact, he was able to extend this result to higher dimensions and gave conditions under which rollings exist in terms of the shapes of the submanifolds, that is, in terms of both intrinsic and extrinsic data.

For reasons unknown to us, this paper seems to have been forgotten over the years. Nomizu's definition of higher dimensional rolling is equivalent to Sharpe's one in Subsection~\ref{subsec:sharpe} and many of his observations have been rediscovered in~\cite[Appendix B]{sharpe97}. Nevertheless, there is no reference to the paper~\cite{nomizu78} in Sharpe's book.

\section{Revival: The two dimensional case and robotics}\label{sec:2dim}

The aim of this section is to put in context the study of the rolling model for the case of two dimensional manifolds, and how they appeared naturally in problems of sub-Riemannian geometry, robotics and geometric control theory. 

\subsection{Rigidity of integral curves in Cartan's distribution}\label{subsec:bryant}

In the celebrated paper~\cite{bryant-hsu}, R. Bryant and L. Hsu studied curves on a manifold $Q$ of dimension $n\geq3$ tangent to a $(2,n)-$distribution ${D}$. The idea was to analyze the space $\Omega_{{D}}(p,q)$ of differentiable curves in $Q$ connecting two points $p,q\in Q$ and being tangent to ${D}$ (called ${D}$-curves by them). The space $\Omega_{{D}}(p,q)$ is endowed with its natural $C^1$ topology. The idea that ${D}$-curves can be ``rigid'' plays a fundamental role in their paper.
\begin{definition}
A ${D}$-curve $\gamma\colon[0,\tau]\to Q$ is {\em rigid} if there is a $C^1$-neighborhood $\mathscr{U}$ of $\gamma$ in $\Omega_{{D}}(\gamma(0),\gamma(\tau))$ so that every $\gamma_1\in\mathscr{U}$ is a reparametrization of $\gamma$. We say that $\gamma$ is {\em locally rigid} if every point of $I=[0,\tau]$ lies in a subinterval $J\subset I$ so that $\gamma$ restricted to $J$ is rigid.
\end{definition}
Their main result goes as follows.
\begin{theorem}[Bryant \& Hsu 1993]\label{th:bryant}
Let ${D}$ be a non-integrable rank $2$ distribution on a manifold $Q$ of dimension $(2+s)\geq3$.  Suppose further that the distribution ${D}_1=[{D},{D}]$ (which has rank $3$) is nowhere integrable. Then there always exist ${D}$-curves that are locally rigid.
\end{theorem}
They give a more precise description of such curves in terms of projections of characteristic curves in a dense subset of the annihilator of ${D}_1$, but stating it precisely would not serve the purposes of this exposition.

For us, the most relevant part of their work is their section on examples, in particular their study of systems of Cartan type and of rolling surfaces. 

Recall that a bracket generating $(2,5)-$distribution is said to be of Cartan type. In other words $D$ is a Cartan distribution if ${D}_1$ has rank $3$ and ${D}_2=[{D}_1,{D}]$ has rank $5$. As a consequence of Theorem~\ref{th:bryant}, they observe that there is exactly a $5$-parameter family of locally rigid ${D}$-curves. In fact they briefly discuss a remarkable geometric behavior occurring in this situation: if $M$ is connected, then any two points of $M$ can be joined by a piecewise smooth ${D}$-curve, whose smooth segments are rigid.

After all these observations, they devote themselves to the analysis of two oriented surfaces $M$ and $\hat M$ endowed with Riemannian metrics rolling one over another without slipping or twisting. Let $F$ and $\hat F$ be the oriented orthonormal frame bundles of $M$ and $\hat M$. Bryant and Hsu considered the ``state space'' manifold $Q=(F\times\hat F)/{\rm SO}(2)$, where ${\rm SO}(2)$ acts diagonally on the Cartesian product. An element in $Q$ is a triple $(x,\hat x;A)$, where $x\in M$, $\hat x\in\hat M$ and $A\colon T_{x}M\to T_{\hat x}\hat M$ is an oriented isometry. Their formulation is as follows. Consider a curve $\gamma\colon[0,\tau]\to Q$ given by $\gamma(t)=(x(t),\hat x(t);A(t))$, then the no-slip condition reads $A(t)(\dot x(t))=\dot{\hat x}(t)$. The no-twist condition requires some more care. Let $e_1,f_1\colon[0,\tau]\to TM$ be a parallel orthonormal frame along the curve $x(t)$ and let
\[
e_2(t)=A(t)(e_1(t)),\quad f_2(t)=A(t)(f_1(t)),
\]
be the orthonormal frame along $\hat x(t)$ obtained via $A$. The rolling has no-twist whenever the moving frame $e_2,f_2$ is also parallel (along $\hat x$).

An important insight for the problem was expressing the no-twist and no-slip conditions in terms of a $(2,5)-$distribution ${D}$ on $Q$. Let $\alpha_1,\alpha_2,\alpha_{21}$ be the canonical $1$-forms of $M$ on $F$ and similarly $\beta_1,\beta_2,\beta_{21}$ for $\hat M$, see~\cite{singer76}. Recall that these forms satisfy the so-called structure equations
\begin{align*}
d\alpha_1&=\alpha_{21}\wedge\alpha_2,&d\beta_1&=\beta_{21}\wedge\beta_2,\\
d\alpha_2&=-\alpha_{21}\wedge\alpha_1,&d\beta_2&=-\beta_{21}\wedge\beta_1,\\
d\alpha_{21}&=\kappa\,\alpha_1\wedge\alpha_2,&d\beta_{21}&=\hat\kappa\,\beta_1\wedge\beta_2,
\end{align*}
where $\kappa$ and $\hat\kappa$ are the Gaussian curvatures of $M$ and $\hat M$ respectively. With all of this, one can consider the distribution $\tilde{{D}}$ on $F\times\hat F$ defined by the Pfaffian equations
\[
\alpha_1-\beta_1=\alpha_2-\beta_2=\alpha_{21}-\beta_{21}=0.
\]
The distribution they were looking for corresponds to the ``push-down'' image of $\tilde{{D}}$ under the submersion $F\times\hat F\to Q$. A smooth curve $\gamma\colon[0,\tau]\to Q$ describes a rolling without slipping or twisting if and only if $\gamma$ is a ${D}$-curve.

A remarkable fact is that the distribution ${D}$ is of Cartan type whenever $\kappa-\hat\kappa\neq0$, which is an open set in $Q$. On this set, the corresponding $5$-parameter family of rigid curves describes the rolling of $\hat M$ against $M$ following geodesics.

\subsection{Non-holonomy in robotics}\label{subsec:marigo}

The traditional modeling of a mechanical system considers configurations (or states) of this mechanical system as points $q$ of a smooth finite-dimensional manifold $M$, and the corresponding velocities $\dot q\in T_qM$ are subject to locally independent constraints
in the Pfaffian form 
\begin{equation}\label{eq:NH0}
A(q)\dot q=0,
\end{equation}
where $A(\cdot)$ is an $m\times n$ matrix of real-valued analytic
functions, where $m<n$. Constraints are said to be  {\it holonomic} if their differential
form given by \eqref{eq:NH0} is integrable. In this case, there exist 
integral submanifolds of dimension $n-m$ that are
invariant. If the constraints are not holonomic at some $q_0\in M$ , then there
will exist an integral submanifold containing $q_0$ of dimension $n-m+k$ with $0<k\leq m$. The integer $k$ is referred to as degree of non-holonomy. If $k=m$, the constraints, and 
by extension the system, are said to be maximally non-holonomic (see~\cite{murray-sastry}). 

There is a more convenient way for control theory to describe the constrained system. If
$G(q)$ denotes a matrix whose columns form a basis for the annihilating distribution of
$A(q)$, then all admissible velocities $\dot q\in A(q)^{\perp}\subset T_qM$ can be
written as linear combinations of the columns of $G(q)$,
\begin{equation}\label{eq:NH1}
\dot q=G(q)w=\sum_{i=1}^{n-m}g_i(q)w_i,
\end{equation}
where $w$ is a vector of {\it quasivelocities} taking values in ${\mathbb R}^{n-m}$.
When quasivelocities can be assigned values at will in time,
functions can be regarded as  {\it control}
inputs of the driftless, linear-in-the-control, nonlinear system defined by \eqref{eq:NH1}.
A physical  {\it actuator} is associated to each control input ,
e.g. a motor for electromechanical systems. The issue of non-holonomy of
the original system, i.e. non-integrability of \eqref{eq:NH0}, can be addressed
by studying the distribution $\Delta$ spanned by the the vector fields $g_i$'s and more precisely the corresponding Lie algebra generated by the $g_i$'s. If the system is maximally non-holonomic (or completely controllable), any two configurations $q$ and $q'$ of its
$n$-dimensional manifold can be connected along the flows of $n-m$ vector fields. From an utilitarian 
engineerÕs viewpoint, the latter definition
may be rephrased as Òan $n$-dimensional non-holonomic system
can be steered at will using less than actuators. This formulation
underscores the appealing fact that devices with reduced
hardware complexity can be used to perform nontrivial tasks, if
non-holonomy is introduced on purpose, and cleverly exploited,
in the device design (see~\cite{murray-sastry}).

Non-holonomy of rolling is particularly relevant to robotic
manipulation, one of the main goals of which is to manipulate an
object grasped by a robot end-effector so as to relocate and re-orient it arbitrarily, the so-called dexterity property. Dexterous
robotic hands developed so far according to an anthropomorphic
paradigm employ far too many joints and actuators (a minimum
of nine) to be a viable industrial solution. Non-holonomy
of rolling can be used to alleviate this limitation. In fact, while
rolling between the surfaces of the manipulated object and that
of fingers has been previously regarded as a complication to be
neglected, or compensated for, some works (see, in particular, \cite{ACL,Bic1,chelouah01,CMP,BM,mar-bic2} and the references therein) tried to exploit rolling for achieving dexterity with simpler
mechanical hardware.

Introducing non-holonomy on purpose in the design of robotic
mechanisms can be regarded as a means of lifting complexity
from hardware to the software and control level of design. In fact,
planning and controlling non-holonomic systems is in general a
considerably more difficult task than for holonomic systems. The
very fact that there are fewer degrees-of-freedom available than
there are configurations implies that standard motion planning
techniques can not be directly adapted to non-holonomic systems.
From the control viewpoint, non-holonomic systems are intrinsically
nonlinear systems, in the sense that they are not exactly
feedback linearizable, nor does their linear approximation retain
the fundamental characteristics of the system, such as controllability (see~\cite{murray-sastry}).

The system of rolling bodies considered here
differs substantially from the class of chained form
systems or differentially flat systems (see Rouchon~\cite{Rou}). Consider, for example, the case such of the {\it plate-ball
system} (i.e. a ball rolling on a plane without slipping or spinning), which is a classical problem in rational mechanics,
brought to the attention of the control community by Brockett
and Dai~\cite{BD}. Montana~\cite{Mon} derived a differential-geometric
model of the rolling constraint between general bodies,
and discussed applications to robotic manipulation. Li and
Canny~\cite{LC} showed that the plate-ball system is controllable,
and that the same holds for two rolling spheres, provided that
their radii are different. 

We close this subsection mentioning the beautiful works of Jurdjevic~\cite{jurd93,jurd95} who studied the problem of
finding the path that minimizes the length of the curve traced
out by the sphere on the fixed plane. It turns out that optimal
paths also minimize the integral of their geodesic curvature, so
that solutions are those of Euler's {\em elastica} problem. For the higher dimensional cases of this problem, see~\cite{jurd08,zimm05}.

\subsection{Orbits and complete answer for controllability}\label{subsec:agrachev}

The point of view adopted by Bryant and Hsu was improved significantly by A. Agrachev and Y. Sachkov in~\cite{agrachev99} employing tools in geometric control theory. 

Two innocent, yet powerful, changes in perspective made the problem more accessible for the application of the orbit theorem of Sussmann~\cite{sussmann73}. These modifications consist of rewriting the state space of the rolling and, most importantly, to prefer the use of vector fields (written in local coordinates) instead of differential forms (written without using coordinates).

Let $M$ and $\widehat M$ be smooth two-dimensional connected oriented Riemannian surfaces. The new version of the state space is given by
\[
Q=\{A\colon T_xM\to T_{\hat x}\widehat M\,|\,x\in M, \hat x\in\widehat M,A\mbox{ an oriented isometry}\}.
\]
It is an easy exercise to see that $Q$ is indeed diffeomorphic to the manifold $M$ introduced in Subsection~\ref{subsec:bryant}. The natural projection $Q\to M\times\widehat M$ is a principal ${\rm SO}(2)$-bundle. As before, a curve $\gamma\colon[0,\tau]\to Q$ describes a rolling motion if there is no slipping, that is, if
\(
A(t)(\dot x(t))=\dot{\hat x}(t)
\)
and there is no twisting (see~\cite{agrachev99})
\[
A(t)\big(\mbox{vector field parallel along }x(t)\big)=\big(\mbox{vector field parallel along }\hat x(t)\big).
\]

Let us now give expressions of the rolling distribution in local coordinates about a point $(x,\hat x;A)\in Q$. Let us consider local orthonormal frames $e_1,e_2$ for $M$ and $\hat e_1,\hat e_2$ for $\hat M$. They define their structure constants $c_1,c_2\in C^\infty(M)$ and $\hat c_1,\hat c_2\in C^\infty(\hat M)$ by the equations $[e_1,e_2]=c_1e_1+c_2e_2$ on $M$ and $[\hat e_1,\hat e_2]=\hat c_1\hat e_1+\hat c_2\hat e_2$ on $\hat M$.

Since $Q$ is a circle bundle over $M\times \hat M$, in the natural trivialization, there a well defined angular direction $\frac{\partial}{\partial\theta}$ and we can identify the isometry $A$ with an angle $\theta$. With these notations, the rolling distribution $D_{\rm R}$ is spanned by the vector fields
\begin{align*}
X_1&=e_1+\cos\theta\,\hat e_1+\sin\theta\,\hat e_2+\big(-c_1+\hat c_1\cos\theta+\hat c_2\sin\theta\big)\frac{\partial}{\partial\theta},\\
X_2&=e_2-\sin\theta\,\hat e_1+\cos\theta\,\hat e_2+\big(-c_2-\hat c_1\sin\theta+\hat c_2\cos\theta\big)\frac{\partial}{\partial\theta}.
\end{align*}

The main controllability theorem for the system of two Riemannian surfaces rolling, as presented in~\cite[Chapter 24]{agrachev04}, is the following.

\begin{theorem}[Agrachev \& Sachkov 1999]
Let ${\mathcal O}={\mathcal O}_{D_{\rm R}}(q)$ be the orbit of the rolling distribution starting at $q\in Q$ and let $\kappa$ and $\hat \kappa$ be the Gaussian curvatures of $M$ and $\hat M$ respectively. Then:
\begin{enumerate}
\item The orbit ${\mathcal O}$ is a imbedded connected submanifold of $Q$ of dimension 2 or 5. More precisely, one has that if $(\kappa-\hat\kappa)|_{\mathcal O}$ is identically zero, then $\dim{\mathcal O}=2$; and if $(\kappa-\hat\kappa)|_{\mathcal O}$ is not identically zero, then $\dim{\mathcal O}=5$.

\item There is an injective correspondence between isometries $\iota\colon M\to\hat M$ and two dimensional orbits of the rolling system. In particular, if the manifolds $M$ and $\hat M$ are isometric, then the rolling model is not completely controllable.

\item If $M$ and $\hat M$ are complete and simply connected, then the correspondence between isometries $\iota\colon M\to\hat M$ and two dimensional orbits of the rolling system is bijective. In particular, the rolling model is completely controllable if and only if the manifolds $M$ and $\hat M$ are not isometric.
\end{enumerate}
\end{theorem}

\section{Re-discovery of the higher dimensional case and interpolation}\label{sec:higerdim}

Here we briefly review the way the higher dimensional problem of rolling manifolds presented to the control theory community and we explain how this was employed in geometric interpolation theory.

\subsection{Sharpe's definition}\label{subsec:sharpe}

Here we present the definition of rolling maps found in the Appendix B of R. W. Sharpe's book~\cite{sharpe97} with some minor modifications.

\begin{definition}\label{def:rollSharpe}
Let $M, \hat{M}$ be $n$-dimensional submanifolds of $\mathbb{R}^{n+\nu}$. Then, a differentiable map $g:[0,\tau] \to {\rm Isom}({\mathbb R}^{n+\nu})$ satisfying the following conditions
\begin{itemize}
\item There is a piecewise smooth curve $x:[0,\tau] \to M$, such that
\begin{itemize}
\item $g(t) x(t) \in \hat{M}$,
\item $T_{g(t)x(t)}\left(g(t) M \right) = T_{g(t)x(t)} \hat{M}$.
\end{itemize}
\item Furthermore, the curve $\hat{x}(t) := g(t) x(t)$ satisfies the following conditions
\begin{itemize}
\item No-slip: $\dot{g}(t) g(t)^{-1} \hat{x}(t) = 0$.
\item No-twist, tangential part:
$d(\dot{g}(t) g(t)^{-1}) T_{\hat{x}(t)}\hat M \subseteq T_{0}(\dot{g}(t) g(t)^{-1}\hat M)^\bot$.
\item No-twist, normal part:
$d(\dot{g}(t) g(t)^{-1}) T_{\hat{x}(t)}\hat M^\bot \subseteq T_{0}(\dot{g}(t) g(t)^{-1}\hat M)$.
\end{itemize}
\end{itemize}
for any $t \in [0,\tau]$ is called a rolling map of $M$ on $\hat{M}$ without slipping or twisting.
\end{definition}

We do not know whether Sharpe was aware of the existence of the paper~\cite{nomizu78} at the time of the publication of his book, but his deduction of the ``correct'' definition rolling maps follows the same structure as Nomizu's. Nevertheless, Sharpe does obtain plenty of extra information. For example he shows that in the imbedded rolling system there is a deep relation with the Levi-Civita connections of the manifolds and the normal connections to the imbeddings. Besides this, he is able to prove precisely that rolling is transitive, that is

\begin{theorem}[Sharpe 1997]
Let $M_0,M_1,M_2\subset{\mathbb R}^{n+\nu}$ be three $n$-dimensional submanifolds, such that they are tangent to each other at a common point $p\in M_0\cap M_1\cap M_2$. Let $\gamma\colon[0,\tau]\to M$ be given such that $\gamma(0)=p$. Assume that $M_1$ rolls on $M_0$ along the curve $\gamma$, with rolling map $g_{1}$, and similarly let $M_2$ roll on $M_1$ along the curve $\hat\gamma=g_{1}\gamma$, with rolling map $g_{2}$. Then $M_2$ rolls on $M_0$ along the curve $\gamma$, with rolling map $g_2g_1$ and with image curve $\tilde\gamma=g_2g_1\gamma=g_2\hat\gamma$.
\end{theorem}

\subsection{Applications to geometric interpolation}\label{subsec:interp}

An interesting application of the rolling system has been in interpolation. The article where this idea appeared for the first time is~\cite{jupp87} for the case of the two dimensional sphere. Afterward it was extended successfully to arbitrary dimensional spheres, Grassmanians and to the special orthogonal groups in~\cite{huper07}. This last application was employed in~\cite{huper06} to study the motion planning of a rotating satellite. Later on in~\cite{huper08} the idea was also shown to work on Stiefel manifolds.

The setting of the interpolation problem seems quite innocent.  Let $x_0,x_1,\dotsc,x_N\in M$ be measurements at times $0=t_0<t_1<\cdots<t_N=\tau$, and consider given initial and final velocities $v\in T_{x_0}M$ and $w\in T_{x_N}M$.
The interpolation problem consists in finding a $C^2$ curve $\gamma\colon[0,\tau]\to M$ satisfying 
\begin{equation*}
\label{int}\gamma(t_i)=x_i,\quad \dot\gamma(0)=v,\,\dot\gamma(\tau)=w,\tag{{\sc Interp}}
\end{equation*}
and $\gamma$ minimizes the functional
\begin{equation*}
\label{func}J(\gamma)=\frac12\int_0^\tau\left\langle\frac{D}{dt}\,\dot\gamma,\frac{D}{dt}\,\dot\gamma\right\rangle\,dt.\tag{{\sc Energy}}
\end{equation*}

Let $\Omega=\{\gamma\in C^2\,|\,\gamma\mbox{ satisfies }\eqref{int}\}$. Then

\begin{theorem}[Crouch \& Silva Leite 1991]\label{ELinterp}
If $\gamma\in\Omega$ minimizes \eqref{func}, then
\begin{equation*}
\frac{D^3}{dt^3}\,\dot\gamma+R\left(\frac{D}{dt}\,\dot\gamma,\dot\gamma\right)\dot\gamma=0,
\end{equation*}
on each $[t_i,t_{i+1}]$, where $R$ is the curvature tensor of $M$.
\end{theorem}

The curves satisfying the differential equation in Theorem~\ref{ELinterp} are called geometric cubic splines, and they are in general quite hard to find. Nevertheless, in the cases described above, the authors were able to find a surprising relation between the rolling dynamics and geometric interpolation. For simplicity of exposition, we only present the relevant results for the case of the $n$ dimensional sphere $S^n$. A first observation that takes place is the following.

\begin{theorem}[Jupp \& Kent 1987, H\"uper \& Silva Leite 2007]
Let $R^\top(t)$ be the rotational part of the rolling map in Definition~\ref{def:rollSharpe}, with rolling curve $\gamma\colon[0,\tau]\to S^n$. For all $t\in[0,\tau]$ and all $j\in{\mathbb N}$,
\[
R^\top(t)\,\frac{D^j}{dt^j}\,\dot\gamma(t)=\gamma_{\rm dev}^{(j+1)}(t),
\]
where $\gamma_{\rm dev}$ is the development of $\gamma$, see Subsection~\ref{subsec:development}.
\end{theorem}

A consequence of the above is the following application to interpolation in $S^n$.

\begin{corollary}[Jupp \& Kent 1987, H\"uper \& Silva Leite 2007]
If the development $t\mapsto\gamma_{\rm dev}(t)$ is an Euclidean cubic spline, then $t\mapsto\gamma(t)$ is a geometric cubic spline on $S^n$ if and only if it is a re-parameterized geodesic.
\end{corollary}

\section{Nowadays: The coordinate-free approach}\label{sec:intrinsic}

The intrinsic definition of the rolling model in higher dimensions was presented for the first time in~\cite{arxiv,norway}. It is clearly motivated by the definition given by Agrachev and Sachkov in~\cite{agrachev99}.

Let $(M,g)$ and $(\hat M,\hat g)$ be two oriented Riemannian manifolds of dimension $n$. The state space of the rolling problem is the manifold
\[
Q=Q(M,\hat M)=\big\{A:T|_x M\to T|_{\hat{x}} \hat{M}\ \big|\ x\in M,\ \hat{x}\in\hat{M}, A \textrm{linear isometry},\ \det(A)>0\big\}.
\]
An absolutely continuous curve $q(t) = (\gamma(t),\hat{\gamma}(t),A(t))$ in $Q$ is a rolling curve if $A(t)X(t)$ is parallel along $\hat{\gamma}(t)$ for every vector field $X(t)$ that is parallel along $\gamma(t)$ (no twist condition) and if $A(t)\dot{\gamma}(t) = \dot{\hat{\gamma}}(t)$ (no slip condition). 

A counting argument shows that $Q$ has dimension $\frac12n(n+3)$. Over this manifold there is an $n$-dimensional distribution $D_{\rm R}$, called the rolling distribution, such that the rolling curves in $Q$ are exactly the integral curves of $D_{\rm R}$. Let us describe this distribution briefly as given in \cite{arxiv}. For a configuration $q=(x,\hat{x};A)\in Q$, and an initial velocity $X\in T|_x M$, we define the \emph{rolling lift} ${\mathscr L}_{\rm R}(X)|_q\in T|_q Q$ as
\begin{align}\label{eq:2.5:3}
{\mathscr L}_{\rm R}(X)|_q=\frac{d}{dt}\big|_0 (P_0^t(\hat{\gamma})\circ A\circ P_t^0(\gamma)),
\end{align}
where $\gamma,\hat{\gamma}$ are any smooth curves in $M,\hat{M}$, respectively, such that $\dot{\gamma}(0)=X$ and $\dot{\hat{\gamma}}(0)=AX$, and $P^b_a(\gamma)$ (resp. $P_a^b(\hat{\gamma})$)
denotes the parallel transport along $\gamma$ from $\gamma(a)$ to $\gamma(b)$
(resp. along $\hat{\gamma}$ from $\hat{\gamma}(a)$ to $\hat{\gamma}(b)$).

\begin{definition} (cf. \cite{arxiv}).
The \emph{rolling distribution} $D_{\rm R}$ on $Q$ is the $n$-dimensional smooth distribution defined, for $q=(x,\hat{x};A)\in Q$, by
$D_{\rm R}|_{q}={\mathscr L}_{\rm R}(T|_x M)|_{q}$.
\end{definition}

An interpretation of the rolling lift ${\mathscr L}_{\rm R}(X)|_q$ of $X\in T|_xM$ at $q=(x,\hat x;A)$ is as follows. Let $\gamma$ be a curve in $M$ such that $\gamma(0)=x$ and $\dot\gamma(0)=X$ then, by the general theory of ordinary differential equations, for short times there is a rolling curve $q(t)$ of $M$ on $\hat M$ satisfying $q(0)=q$. The rolling lift is precisely $\dot q(0)$.

\subsection{The controllability problem}\label{subsec:control}

The orbit ${\mathcal O}_{D_{\rm R}}(q)$ of the rolling system described above passing through $q\in Q$ consists of all the states $\tilde q$ that can be connected to $q$ via a rolling curve. The (complete) controllability problem asks for conditions on the geometry of $M$ and $\hat M$ such that ${\mathcal O}_{D_{\rm R}}(q)=Q$. One way of addressing this problem is via Sussmann's orbit theorem, that is, by showing that all the Lie brackets of the vector fields steering the dynamics have to span the tangent bundle of the state space. For the rolling model, this Lie brackets are expressed in terms of the curvature tensors $R$ and $\hat{R}$ associated to the Riemannian metrics $g$ on $M$ and  $\hat{g}$ on $\hat{M}$ respectively, together with the covariant derivatives of $R$ and $\hat{R}$. It seems therefore impossible to solve for general dimension $n$ the controllability issue on the sole knowledge of the Lie algebraic structure of $D_{\rm R}$, except for low dimensions. Indeed, in the case for instance where $(\hat{M},\hat{g})$ is the $n$-dimensional Euclidean space, it would amount to determine ${\rm Hol}(\nabla^g)$, the holonomy group of the Levi-Civita connection $\nabla^g$ associated to $g$, with the only knowledge of its curvature tensor and its covariant derivatives. Instead, the latter issue can be successfully addressed by resorting on group theoretic and algebraic arguments, see~\cite{arxiv}. For specific examples, using extra knowledge of the problem at hand, see~\cite{norway,jurd08,zimm05}.

In general, one can define a notion of curvature especially adapted to the rolling model, see~\cite{CK}. For $q=(x,\hat x;A)\in Q$, the \emph{rolling curvature} is the linear map
\[
{\rm Rol}_q\colon\bigwedge^2T|_xM\to T^*|_xM\otimes T|_{\hat x}\hat M;
\quad {\rm Rol}_q(X\wedge Y):=AR(X,Y)-\hat R(X,Y)A.
\]
This map permits to give a first sufficient condition for the rolling model to be controllable, see~\cite{arxiv,grong}. 

\begin{theorem}[Chitour \& Kokkonen 2011, Grong 2012]
If the \emph{rolling curvature} is an isomorphism 
for every $q\in Q$, then the rolling problem is completely controllable.
\end{theorem}

The above condition is very hard to deduce directly from conditions on the geometry of $M$ and $\hat M$. It is therefore necessary to reduce the problem to a simpler one. One possible way to do this is to give some extra structure to the manifold $\hat M$. In this vein, it was possible to give controllability conditions ``without Lie brackets'' for the case in which $(\hat M,\hat g)=({\mathbb F}^n_c,{\bf g}^n_c)$ is the space form of constant sectional curvature $c$, see~\cite{kobayashi63,sakai91}. To state these, let us first introduce some terminology.

\begin{definition}
Consider the vector bundle $\pi_{TM\oplus{\mathbb R}}:TM\oplus{\mathbb R}\to M$. The rolling connection $\nabla^c$ is the vector bundle connection on $\pi_{TM\oplus{\mathbb R}}$ defined by
\begin{align}\label{eq:nabla_rol_explicit}
\nabla^{c}_X (Y,s)=\Big(\nabla_X Y+s(x)X,X(s)-cg\big(Y|_x,X)\Big),
\end{align}
for every $x\in M$, $X\in T|_x M$, $(Y,s)\in{\rm VF}(M)\times C^\infty(M)$; where we have canonically identified the space of smooth sections $\Gamma(\pi_{TM\oplus{\mathbb R}})$ of $\pi_{TM\oplus{\mathbb R}}$ with ${\rm VF}(M)\times C^\infty(M)$.
\end{definition}

When $c\neq0$, the connection $\nabla^{c}$ is a metric connection with respect to the fiber inner product $h_c$ on $TM\oplus{\mathbb R}$ defined by
\[
h_c((X,r),(Y,s))=g(X,Y)+c^{-1}rs,
\]
where $X,Y\in T|_x M$, $r,s\in{\mathbb R}$. Its holonomy group is denoted by ${\mathcal H}^{c}(M)$. In this language, we have the following result, see~\cite{CK}.

\begin{theorem}[Chitour \& Kokkonen 2012]\label{holcontr}
Let $(M,g)$ be a complete, oriented and simply connected Riemannian manifold. The rolling problem of $M$ rolling against ${\mathbb F}^n_c$ is completely controllable if and only if
\[
{\mathcal H}^{c}(M)=\left\{\begin{array}{lc}
{\rm SO}(n+1),&c>0;\\
{\rm SE}(n),&c=0;\\
{\rm SO}_0(n,1),&c<0.
\end{array}\right.
\]
Here the Lie group ${\rm SO}_0(n,1)$ represents the identity component of the group ${\rm O}(n,1)$ of linear transformations that preserve the quadratic form $F_{n,1}(x_1,\dotsc,x_{n+1})=x_1^2+\cdots+x_n^2-x_{n+1}^2$.
\end{theorem}

Wanting to fully understand these cases, it is important to remark some structure theorems encoded in Theorem~\ref{holcontr}. Observe that up to rescaling, it is sufficient to study when $c=0$, $1$ and $-1$. In the Euclidean situation, i.e. $c=0$, the condition ${\mathcal H}^{0}(M)={\rm SE}(n)$ is equivalent to that $M$ has full Riemannian holonomy. In the case $c=1$, if the action of ${\mathcal H}^{1}(M)$ on the unit sphere is not transitive, then $(M, g)$ is the unit sphere. As a consequence, it holds that, for $n\geq 16$ and even, the rolling system $Q=Q(M,S^{n-1})$ is completely controllable if and only if $(M,g)$ is not isometric to the unit sphere. Both theses cases were analyzed in~\cite{CK}, and the remaining cases are currently under investigation. The hyperbolic case presented a more difficult challenge, see~\cite{CGK1}.

\begin{theorem}[Chitour, Godoy \& Kokkonen 2012]
Let $(M,g)$ be a complete, oriented and simply connected Riemannian $n$-manifold rolling onto the space form $({\mathbb H}^n,{\bf g}^n_{-1})$ of curvature $-1$. Then the associated rolling problem is completely controllable if and only if $(M,g)$ is not isometric to a warped product of the form
\begin{description}
\item[{\rm (WP1)}] $({\mathbb R}\times M_1,ds^2\oplus_{e^{cs}} g_1)$, or

\item[{\rm (WP2)}] $({\mathbb H}^k\times M_1,{\bf g}^{k}_{-1}\oplus_{\cosh(\sqrt{-c}\,d)} g_1)$, where $1\leq k\leq n$ and for each $x\in {\mathbb H}^k$ , $d(x)$ is the distance between $x$ and an arbitrary fixed point $x_0\in{\mathbb H}^k$. 
\end{description}
In both situations, $(M_1,g_1)$ is some complete simply connected Riemannian manifold. As usual, the term $ds^2$ represents the usual Riemannian metric on ${\mathbb R}$.
\end{theorem}

\subsection{Symmetries of the rolling problem}\label{subsec:symmetries}

The idea developed in Subsection~\ref{subsec:control} of setting $\hat M$ to be a space form has a beautiful geometric consequence on the bundle structure of the natural projection $\pi_{Q,M}\colon Q=Q(M,{\mathbb F}^n_c)\to M$. Let us explain what this is.

In general, it is not clear if there is a $G$-principal bundle structure on $Q$ making $D_{\rm R}$ a $G$-principal bundle connection. This is indeed the case if the manifolds are of dimension 2, in which case the projection $Q\to M\times\hat M$ is a principal ${\rm SO}(2)$ bundle with $D_{\rm R}$ as its connection. Nevertheless, for higher dimensions the projection $Q\to M\times\hat M$ does not satisfy the above. The main reason is that ${\rm SO}(n)$ is abelian only for $n=2$, thus we need to make the problem simpler.

For $c\neq0$, let $G_c(n)$ be the identity component of the Lie group of linear maps that leave invariant the bilinear form $\langle x,y\rangle^n_{c}:=\sum_{i=1}^n x_i y_i+c^{-1}x_{n+1}y_{n+1}$, for $x=(x_1,\dots,x_{n+1}),y=(y_1,\dots,y_{n+1})\in{\mathbb R}^{n+1}$. Observe that $G_1(n)={\rm SO}(n+1)$ and $G_{-1}(n)={\rm SO}_0(n,1)$. For $c=0$, we set $G_{0}(n)={\rm SE}(n)$. Recall that, with this notation, the identity component of the isometry group of $({\mathbb F}^n_c,{\bf g}^n_{c})$ is equal to $G_c(n)$ for all $c\in{\mathbb R}$ (cf. \cite{kobayashi63}).

The fundamental result concerning rolling against a space form lies in the fact that there is a $G_c(n)$-principal bundle structure for the state space compatible with the distribution $D_{\rm R}$, i.e. $D_{\rm R}$ is a $G_c(n)$-principal bundle connection, see~\cite{CK}. The precise result follows.

\begin{theorem}[Chitour \& Kokkonen 2012]
Let $Q=Q(M,{\mathbb F}^n_c)$ be the state space of rolling $M$ against the space form ${\mathbb F}^n_c$. Then we have:
\begin{itemize}
\item[(i)]
The projection $\pi_{Q,M}:Q\to M$ is a principal $G_c(n)$-bundle with a left action $\mu:G_c(n)\times Q\to Q$ defined for every $q=(x,\hat{x};A)$ by
\[
\mu((\hat{y},C),q)=(x,C\hat{x}+\hat{y};C\circ A), \textrm{if $c=0$}, 
\mu(B,q)=(x,B\hat{x};B\circ A), \textrm{if $c\neq 0$}. 
\]
Moreover, the action $\mu$ preserves the distribution $D_{\rm R}$, i.e., for any $q\in Q$ and $B\in G_c(n)$, $(\mu_B)_*D_{\rm R}|_q=D_{\rm R}|_{\mu(B,q)}$, where $\mu_B:Q\to Q$; $q\mapsto \mu(B,q)$.

\item[(ii)] For any given $q=(x,\hat{x};A)\in Q$, there is a unique subgroup ${\mathcal H}^{c}_q$ of $G_c(n)$, called the holonomy group of $D_{\rm R}$ at $q\in Q$, such that 
\[
\mu({\mathcal H}^{c}_q\times\{q\})={\mathcal{O}}_{D_{\rm R}}(q)\cap \pi_{Q,M}^{-1}(x).
\]
Also, if $q'=(x,\hat{x}';A')\in Q$ is in the same $\pi_{Q,M}$-fiber as $q$, then ${\mathcal{H}}^c_q$ and ${\mathcal{H}}^c_{q'}$ are conjugate in $G_c(n)$ and all conjugacy classes of ${\mathcal{H}}^c_q$ in $G_c(n)$ are of the form~${\mathcal{H}}^c_{q'}$.
\end{itemize}
\end{theorem}

A natural question to ask is whether a  converse of the theorem above holds, in other words, does the existence of a $G$-principal bundle structure on $Q$ such that $D_{\rm R}$ is a connection imply that $\hat M$ must have constant sectional curvature? The answer is generically yes, but we need to introduce some more terminology.

Recall that in Subsection~\ref{subsec:1/3} we defined the Lie algebra of symmetries ${\rm Sym}(D)$ of a distribution $D$ on a manifold $\tilde M$ as the set of vector fields $X\in{\rm VF}(\tilde M)$ that satisfy $[X,D]\subseteq D$. For the case of the rolling distribution, we will focus our attention in the symmetries of the rolling distribution that are annihilated by the projection $\pi_{Q,M}\colon Q\to M$, that is, in the Lie algebra
\[
{\rm Sym}_0(D_{\rm R}):=\{S\in{\rm Sym}(D_{\rm R})\ |\ (\pi_{Q,M})_*S=0\}.
\]
With this at hand, the mentioned converse takes the following form, see~\cite{CGK2}.

\begin{theorem}[Chitour, Godoy \& Kokkonen 2012]
If there is an open dense set $O\subset Q$
such that $R|_x:\bigwedge^2 T|_x M\to \bigwedge^2 T|_x M$ is invertible on $\pi_{Q,M}(O)$
and $\widetilde{\rm Rol}$ is invertible on $O$,
then, up to an isomorphism of Lie-algebras,
\[
{\rm Sym}_0(D_{\rm R})=\mathrm{Iso}(\hat{M},\hat{g})
\]
and therefore all the elements of ${\rm Sym}_0(D_{\rm R})$ are induced by Killing fields of $(\hat{M},\hat{g})$.

In particular, under the above assumptions, if
there is a principal bundle structure on $\pi_{Q,M}:Q\to M$
that renders $D_{\rm R}$ to a principal bundle connection,
then $(\hat{M},\hat{g})$ is a space of constant curvature.
\end{theorem}

\subsection{Generalizations and perspectives}\label{subsec:generalization}

Two natural questions to ask concern the extension of the rolling system to the situation in which the manifolds involved have different dimension and to to extend the classification result in Subsection~\ref{subsec:agrachev} to other cases. For the first question, one needs to consider curves of isometric injections instead of isometries. This change introduces many difficulties in understanding the controllability problem, and in fact many tools that work well in the classical situation can not be generalized. The second question has a satisfactory answer for the three dimensional case, see~\cite{arxiv}. There it is shown that the orbits can have dimensions 3, 6, 7, 8 and 9.

A question that has been in our minds for a while is to actually compare the manifolds via the rolling model. This idea of comparison is naively evident in the rolling curvature tensor: one is actually subtracting the Riemannian curvatures of the manifolds. In fact, rolling should provide a framework for the isometric characterization of manifolds by using curvature tensor spectrum information (as in Osserman-type conditions for instance, cf. \cite{Gilkey}).

Finally, we have noticed that the problem of rolling manifolds can be generalized far beyond than allowing arbitrary connections, as in~\cite[Section 7]{grong}, or to pseudo-Riemannian manifolds, as in~\cite{MS}. This extension consists of rolling so-called Cartan geometries, see~\cite{sharpe97}, and it includes as particular cases both of the situations mentioned above, together with the problem of rolling manifolds of different dimensions, see~\cite{CGK3}. The main idea behind this is that Cartan geometries are the most general framework for a notion of development to exists, which underlies the very definition of the rolling dynamics. So far this generalized model has resisted a thorough study of controllability.

\end{document}